\documentclass[preprint, 12pt, 3p]{article}

\usepackage{amsmath, amssymb, amsfonts, amsthm}
\RequirePackage[numbers]{natbib}
\RequirePackage[colorlinks,citecolor=blue,urlcolor=blue]{hyperref}

\newtheorem{theorem}{Theorem}

\newtheorem{definition}{Definition}

\newtheorem{example}{Example}

\newtheorem{Assumption}{Assumption}
\newtheorem{Remark}{Remark}

\newcommand{\slim} {\mathop{\rm lim\,sup\,}}
\newcommand{\ilim} {\mathop{\rm lim\,inf\,}}

\def\S{\mathbb{S}}

\def\Y{\mathbb{Y}}
\def\h{\mathbf{I}}

\def\A{\mathbb{A}}

\def\R{\mathbb{R}}
\def\P{\mathbb{P}}
\def\F{\mathbb{F}}

\def\W{\mathbb{W}}

\def\c{\bar{c}}
\def\B{\mathcal{B}}
\def\oo{\mathcal{O}}

\def\fff{\mathtt{F}}

\def\Psii{{\Xi_{\int}}}

\begin{document}


\title{Equivalent Conditions for Weak Continuity of Nonlinear Filters}

\maketitle

%

\begin{center}
Eugene~A.~Feinberg \footnote{Department of Applied Mathematics and
Statistics,
 Stony Brook University,
Stony Brook, NY 11794-3600, USA, eugene.feinberg@sunysb.edu}\ and
Pavlo~O.~Kasyanov\footnote{Institute for Applied System Analysis,
National Technical University of Ukraine ``Kyiv Polytechnic
Institute'', Peremogy ave., 37, build, 35, 03056, Kyiv, Ukraine,\
kasyanov@i.ua.}

\end{center}

\begin{abstract}
This paper studies weak continuity of nonlinear filters. It is well-known that Borel measurability of transition probabilities for problems with incomplete state observations is preserved when the original discrete-time process is replaced with the process whose states are belief probabilities.  It is also known that the similar preservation may not hold for weak continuity of transition probabilities.  In this paper we show that the sufficient condition for weak continuity of transition probabilities for beliefs introduced by Kara, Saldi, and Y\"uksel (2019) is a necessary and sufficient condition for semi-uniform Feller continuity of transition probabilities.  The property of semi-uniform Feller continuity was introduced recently by Feinberg, Kasyanov, and Zgurovsky (2022), and  the original transition probability for a Markov decision processes with incomplete information has this property if and only if the transition probability of the process, whose state is a pair consisting of the belief probability and observation, also has this property. Thus, this property implies weak continuity of nonlinear filters.  This paper also reviews several  necessary and sufficient conditions for semi-uniform Feller continuity.
\end{abstract}



\section{Introduction}  \label{s1} As was understood long ago in \cite{Ao,As,Dy,Shi64,Shi}, the main general method for studying problems with incomplete information is their reduction to problems with belief states or, in other words, posterior distributions of the states.  This is true  for problems with Borel state, observation, and action spaces \cite{Rh,Yu}. However, an important property for stochastic optimization is weak continuity of transition probabilities, and this property may not hold for the process with belief states even if it holds for the original process \cite[Example 4.1]{FKZ}.

This paper studies problems with a  hidden state set $\W,$ a set of observations $\Y,$ and a set of decisions (or controls) $\A.$ These sets are Borel subsets of Polish (complete separable metric) spaces.  We consider four models: a \emph{Markov Decision Process with Incomplete Information (MDPII)}, \emph{Platzman's model}, and two models of \emph{Partially Observable Markov Decision Processes (POMDPs)}:  \emph{${\rm POMDP}_1$} and  \emph{${\rm POMDP}_2.$ } An MDPII, also known under several other names, is probably the oldest model. This model and its versions are described in many references including monographs  \cite{BR,BS,DY} and mentioned above references  \cite{Ao,As,Dy,Rh,Shi64,Shi,Yu}.

The dynamics of an MDPII is defined by transition probabilities\\ $P(dw_{t+1}dy_{t+1}|w_t,y_t,a_t),$ where  $w_t\in\W$ is the hidden state, $y_t\in\Y$ is the observation, and  $a_t\in\A$ is the selected control at the time epoch $t=0,1,\ldots.$
Platzman'a model is an MPDII, for which transition probabilities do not depend on observations, that is, $P(dw_{t+1}dy_{t+1}|w_t,y_t,a_t)=P(dw_{t+1}dy_{t+1}|w_t,a_t).$ This model was introduced in \cite{Plat}, where it was observed that  two different models of POMDPs had been studied in the literature.  These models were called ${\rm POMDP}_1$ and  ${\rm POMDP}_2$ in \cite{FKZSIAM}.

${\rm POMDP}_1$ is a Platzman's model with the transition probability \\ $P(dw_{t+1}dy_{t+1}|w_t,a_t)=P_1(dw_{t+1}|w_t,a_t)
Q_1(dy_{t+1}|w_t,a_t),$ $t=0,1,\ldots,$ where $P_1$ is the transition probability for hidden states, and $Q_1$ is the observation probability.   ${\rm POMDP}_2$ is  Platzman's model with the transition probability $P(dw_{t+1}dy_{t+1}|w_t,a_t)=P_2(dw_{t+1}|w_t,a_t)
Q_2(dy_{t+1}|a_t,w_{t+1}),$ $t=0,1,\ldots,$ where $P_2$ is the transition probability for hidden states, and $Q_2$ is the observation probability. A ${\rm POMDP}_1$ is mostly used in operations research, and ${\rm POMDP}_2$ is used both in operations research and electrical engineering, and this model describes nonlinear Kalman filters; see \cite{FKZ, FKZSIAM, HL, Plat} for details. For infinite-state problems, most of the results on continuity of transition probabilities for beliefs are currently known for ${\rm POMDP}_2$  \cite{Steklov, FKZ, FKZSIAM, HL, Saldi, RS}.

For ${\rm POMDP}_2$ sufficient conditions for weak continuity of transition probabilities for beliefs are provided in monographs \cite[p. 92]{HL} and \cite[ Chapter 2]{RS}. They both assume weak continuity of transition probabilities $P_2$ and continuity in total variations of the observation probabilities $Q_2.$ They assumed other additional conditions.  In \cite{FKZ} it was shown that  weak continuity of transition probabilities $P_2$ and continuity in total variations of the observation probabilities $Q_2$ are sufficient for weak continuity of transition probabilities for beliefs.  This was done by using the uniform Fatou lemma~\cite{UFL} and Assumption~\ref{Ass:H} below on continuity properties of transition probabilities. Another proof of this fact was provided in \cite{Saldi}, where it was also provided another sufficient condition  for weak continuity of transition probabilities for belief states; see assumption (iii) in Section~\ref{s4}.  In addition, a more general assumption (see Assumption~\ref{Ass:M} in Section~\ref{s2} and Assumption (M) in Section~\ref{s4}) is provided in \cite{Saldi} as an apparently simpler alternative to Assumption~\ref{Ass:H}.

Sufficient conditions for weak continuity of transition probabilities for completely observable Markov Decision Processes corresponding to MDPIIs were studied in \cite{Steklov,FKZ,FKZSIAM}. Such completely observable models are called Markov Decision Processes with Complete Information (MDPCIs).  A state of an MDPCI is a pair $(z_t,y_t),$ where $z_t$ is the belief probability (posterior probability distribution of the state) and $y_t$ is the observation at epach $t=0,1,\ldots.$ A state of belief-MDPs, which can be constructed for Platzman's models and POMDPs, is the posterior probability distribution of the state $z_t,$ $t=0,1,\ldots.$  States $(z_t,y_t)$ can be also considered for models with complete information for Platzman's models and for POMDPs  since by definitions MDPIIs are more general models than Platzman's models and POMDPs. In this paper the transition probability for a completely observable model with states   $(z_t,y_t)$ is denoted by $q(dz_{t+1}dy_{t+1}|z_t,y_t,a_t),$ and its marginal distribution is  $\tilde{q}(dz_{t+1}|z_t,y_t,a_t):=q(dz_{t+1},\Y|z_t,y_t,a_t).$  Therefore, weak continuity of $q$ implies weak continuity of $\tilde{q}.$  For Platzman's models and POMDPs these transition probabilities do not depend on observations $y_t,$  and $\tilde{q}$ is the transition probability of the belief-MDP.


Continuity of belief probabilities for MDPCIs were studied in \cite{Steklov, FKZ}, and recently MDPCIs with semi-uniform Feller transition probabilities and their applications to  Platzman's models and POMDPs were investigated in \cite{FKZSIAM}. The notion of semi-uniform Feller transition probabilities was introduced in \cite{FKZJOT}. This property is stronger than weak continuity.  
This property provides the straightforward way to prove weak continuity of the transition probability $\hat{q}$ for belief-MDPs for some problems..  As shown in \cite{FKZSIAM}, the transition probability $q$  for beliefs is semi-uniform Feller  if and only the original transition probability $P$ is semi-uniform Feller; see Theorem~\ref{th:extra} below. Semi-uniform Feller continuity of $q$ implies weak continuity of $q.$  Weak continuity of $q$ implies weak continuity of $\hat{q}.$
 In addition, in view of Theorem~\ref{th:extra} below, semi-uniform Feller continuity of the kernel $P$  is equivalent to semi-uniform Feller continuity of the kernel $R,$ which is an integrated version of the kernel $P$ defined in \eqref{3.3} for MDPIIs and in \eqref{3.3ef} for Platzman's models and POMDPs.

Therefore, a natural research direction is to identify necessary and sufficient conditions for semi-uniform Feller continuity of a transition kernel. 
Two such conditions, were introduced in \cite{FKZJOT}. The first necessary and sufficient condition is  Assumptions~\ref{AssKern} stated below.  The second one is Assumption~\ref{Ass:H} taken together with continuity of the margin kernel; see Theorem~\ref{th:CBRmain}.   These two conditions are based on sufficient conditions for weak continuity of $\tilde{q}$ for  \emph{${\rm POMDP}_2$}  introduced in \cite{Steklov, FKZ} before semi-uniform Feller continuity was defined in \cite{FKZJOT}.


 This paper introduces the necessary and sufficient Assumption~\ref{Ass:M} based on assumption~(M) introduced in  Kara et al.~\cite{Saldi} as a sufficient condition of weak continuity of $\tilde{q}$ for  \emph{${\rm POMDP}_2.$} 
 As we discussed above, in order to prove weak continuity of the transition kernels $q$ and $\tilde{q},$ it is sufficient to verify semi-uniform continuity of $P.$ This can be done by verifying one of these assumptions for the transition kernel $P.$ 


Section~\ref{s2} of this paper describes properties of semi-uniform Feller kernels.  Theorem~\ref{th:CBRLmain} is the main result of this paper. Section~\ref{s3} describes results on semi-uniform Feller continuity of transition probabilities $q$ for MDPCIs, and Section~\ref{s4} describes sufficient conditions for weak continuity of transition probabilities $\hat{q}$ for belief-MDPs corresponding to Platzman's models and POMDPs.

\section{Semi-Uniform Feller Stochastic Kernels}\label{s2}
For a separable metric space $\S=(\S,\rho_\S),$ where $\rho_\S$ is a metric, let $\tau(\S)$ be the topology of $\S$ (the family of all open subsets of $\S$), and let ${\mathcal B}(\S)$ be its Borel
$\sigma$-field, that is, the $\sigma$-field generated by all open subsets of the
metric space $\S$.
For a subset $S$ of $\S$ let $\bar{S}$ denote the \textit{closure of} $S,$ and $S^o$ is the \textit{interior of} $S.$ Then $S^o$ is open, $\bar{S}$ is closed, and $S^o\subset S\subset \bar{S}.$ Let $\partial S:=\bar{S}\setminus S^o$ denote the \textit{boundary of} $S.$ 

We denote by $\P(\S)$ the \textit{set of probability
measures} on $(\S,{\mathcal B}(\S)).$ A sequence of probability
measures $\{\mu^{(n)}\}_{n=1,2,\ldots}$ from $\P(\S)$
\textit{converges weakly} to $\mu\in\P(\S)$ if for any
bounded continuous function $f$ on $\S$
\begin{equation}\label{eq:wcEF}\int_\S f(s)\mu^{(n)}(ds)\to \int_\S f(s)\mu(ds) \qquad {\rm as \quad
}n\to\infty.
\end{equation}
This definition of weak convergence also applies to a sequence of measures converging to a finite measure $\mu,$  that is, $\mu(\S)<\infty.$
A sequence of probability measures $\{\mu^{(n)}\}_{n=1,2,\ldots}$ from $\P(\S)$ \textit{converges in  total  variation} to $\mu\in\P(\S)$ if
\begin{equation}\label{eq:Kara1}
\begin{aligned}
\sup_{C\in \B(\S)}|\mu^{(n)}(C)-\mu(C)|\to 0\ {\rm as} \ n\to \infty.
\end{aligned}
\end{equation}

Note that $\P(\S)$ is a separable metrizable topological space with respect to the topology of weak convergence for probability measures when $\S$ is
a separable metric space \cite[Chapter~II]{Part}, and there are several ways to introduce a metric on $\P(\S)$ generating this topology.

For a Borel subset $S$ of a metric space $(\S,\rho_\S)$, where $\rho_\S$
is a metric, we always consider the
 metric space $(S,\rho_S),$ where $\rho_S:=\rho_\S\big|_{S\times S}.$  A subset $B$ of $S$ is called open
(closed) in $S$ if $B$ is open (closed) in $(S,\rho_{\color{black}S})$. Of course, if $S=\S$, we omit
``in $\S$''. Observe that, in general, an open (closed) set
in $S$ may not be open (closed). For $S\in\B(\S)$ we denote by $\B(S)$ the Borel $\sigma$-field on
$(S,\rho_S).$  Observe that $\B(S)=\{S\cap B:B\in\B(\S)\}.$
For metric spaces $\S_1$ and $\S_2$, a (Borel-measurable) \textit{stochastic
kernel} $\Psi(ds_1|s_2)$ on $\S_1$ given $\S_2$ is a mapping $\Psi(\,\cdot\,|\,\cdot\,):\B(\S_1)\times \S_2\mapsto [0,1]$ such that $\Psi(\,\cdot\,|s_2)$ is a
probability measure on $\S_1$ for any $s_2\in \S_2$, and $\Psi(B|\,\cdot\,)$ is a Borel-measurable function on $\S_2$ for any Borel set $B\in\B(\S_1)$.  Another name for a stochastic kernel is a transition probability. A
stochastic kernel $\Psi(ds_1|s_2)$ on $\S_1$ given $\S_2$ defines a Borel measurable mapping $s_2\mapsto \Psi(\,\cdot\,|s_2)$ of $\S_2$ to the metric space
$\P(\S_1)$ endowed with the topology of weak convergence.
A stochastic kernel
$\Psi(ds_1|s_2)$ on $\S_1$ given $\S_2$ is called
\textit{weakly continuous (continuous in
total variation)}, if $\Psi(\,\cdot\,|s^{(n)})$ converges weakly (in
total  variation) to $\Psi(\,\cdot\,|s)$ whenever $s^{(n)}$ converges to $s$
in $\S_2$. 

\begin{definition}\label{defi:equi} (\cite{FKL2})
A set $\mathtt{F}$ of real-valued functions on a metric space $\S$ is called
\begin{itemize}
\item[{\rm(i)}]
\textit{lower semi-equicontinuous at a point} $s\in \S$ if $\ilim\limits_{s'\to s}\inf\limits_{f\in\mathtt{F}}(f(s')-f(s))\ge0;$
\item[{\rm(ii)}]
\textit{upper semi-equicontinuous at a point} $s\in \S$ if the set $\{-f\,:\,f\in\mathtt{F}\}$ is lower semi-equicontinuous at $s\in \S;$
\item[{\rm(iii)}] \textit{equicontinuous at a point $s\in\S$}, if $\mathtt{F}$ is both lower and upper semi-equicontinuous at $s\in\S,$ that is, $\mathop{\lim}\limits_{s'\to s} \mathop{\sup}\limits_{f\in\mathtt{F}} |f(s')-f(s)|=0;$
\item[{\rm(iv)}]
\textit{lower / upper semi-equicontinuous (equicontinuous respectively)} (\textit{on $\S$}) if it is lower / upper semi-equicontinuous (equicontinuous respectively) at all $s \in \S;$
\item[{\rm(v)}] \textit{uniformly bounded (on $\S$)}, if there exists a constant $L<+\infty $ such that $ |f(s)|\le L$ for all $s\in\S$ and  for
all $f\in \mathtt{F}.$
\end{itemize}
\end{definition}

Let $\S_1,$ $\S_2,$ and $\S_3$ be Borel subsets of Polish spaces, and $\Psi$ on $\S_1\times\S_2$ given $\S_3$ be a stochastic kernel. For $A\in\B(\S_1),$
$B\in\B(\S_2),$ and $s_3\in\S_3,$ let
\begin{equation}\label{eq:marg_new}
\Psi(A,B|s_3):=\Psi(A\times B|s_3).
\end{equation}
In particular, we consider \textit{marginal} stochastic kernels
$\Psi(\S_1,\,\cdot\,|\,\cdot\,)$ on $\S_2$ given $\S_3$ and $\Psi(\,\cdot\,,\S_2|\,\cdot\,)$ on $\S_1$ given $\S_3.$
%
\begin{definition}\label{defi:unifFP} (\cite{FKZJOT})
A stochastic kernel $\Psi$ on $\S_1\times\S_2$ given $\S_3$ is \emph{semi-uniform Feller} if, for each sequence $\{s_3^{(n)}\}_{n=1,2,\ldots}\subset\S_3$ that converges to $s_3$ in $\S_3$ and for each bounded continuous function $f$ on $\S_1,$
\begin{equation}\label{eq:equivWTV3}
\lim_{n\to\infty} \sup_{B\in \B(\S_2)} \left| \int_{\S_1} f(s_1) \Psi(ds_1,B|s_3^{(n)})-\int_{\S_1} f(s_1) \Psi(ds_1,B|s_3)\right|= 0.
\end{equation}
\end{definition}
A semi-uniform Feller stochastic kernel $\Psi$ on $\S_1\times \S_2$ given $\S_3$ is weakly continuous \cite{FKZJOT,FKZSIAM}. 
%
We recall that the marginal measure $\Psi(ds_1,B|s_3),$ $s_3\in\S_3,$ is defined in \eqref{eq:marg_new}. As follows from \eqref{eq:equivWTV3},  if $\Psi$ is a semi-uniform Feller stochastic kernel on $\S_1\times\S_2$ given $\S_3,$ then for each $B\in\mathcal{B}(\S_2)$ the kernel $\Psi(ds_1,B|s_3)$ on $\S_1$ given $\S_3$ is weakly continuous, that is, if $s^{(n)}_3\to s_3$ as $n\to\infty,$ where $s^{(n)}_3,s_3\in\S_3$ for $n=1,2,\ldots,$ then sequence of substochastic measures $\{\Psi(ds_1,B|s^{(n)}_3)\}_{n=1}^\infty$ converges weakly to $\Psi(ds_1,B|s_3).$

For each set $A\in\B(\S_1)$ consider the set of functions
\begin{equation}\label{eq:familyoffunctions}
\fff^\Psi_A=\{  s_3\mapsto \Psi(A\times B |s_3):\, B\in \B(\S_2)\}
\end{equation}
mapping $\S_3$  into $[0,1].$ Consider the following type of continuity for stochastic kernels on $\S_1\times\S_2$ given $\S_3.$
\begin{definition}\label{defi:wtv}
A stochastic kernel $\Psi$ on $\S_1\times\S_2$ given $\S_3$ is called \textit{WTV-continuous}, if for each $\oo \in\tau (\S_1)$
the set of
functions $\fff^\Psi_\oo$ is lower semi-equicontinuous on $\S_3.$
\end{definition}
 Definition~\ref{defi:equi}{\color{black}(i)} directly implies that the stochastic kernel $\Psi$ on $\S_1\times\S_2$ given $\S_3$ is WTV-continuous if and only if for each $\oo \in\tau (\S_1)$
\begin{equation}\label{eq:equivWTV0}
\ilim_{n\to\infty} \inf_{B\in \B(\S_2)\setminus\{\emptyset\}} \left( \Psi(\oo \times B|s_3^{(n)})-\Psi(\oo \times B|s_3)\right)\ge 0,
\end{equation}
whenever $s_3^{(n)}$ converges to $s_3$ in $\S_3.$  
{\color{black} Since} $\emptyset\in\B(\S_2),$ \eqref{eq:equivWTV0} holds if and only if
\begin{equation}\label{eq:wtvsc}
\lim_{n\to\infty} \inf_{B\in \B(\S_2)} \left( \Psi(\oo \times B|s_3^{(n)})-\Psi(\oo \times B|s_3)\right)= 0.
\end{equation}

The following theorem provides necessary and sufficient conditions for semi-uniform Feller continuity of stochastic kernels; see the relevant facts for weak continuity in \cite{Part,sha}.

\begin{theorem}\label{th:equivWTV} (\cite{FKZJOT})
For a stochastic kernel $\Psi$ on $\S_1\times\S_2$ given $\S_3,$  the following conditions are equivalent:
\begin{itemize}
\item[{\rm(a)}] the stochastic kernel $\Psi$ on $\S_1\times\S_2$ given $\S_3$ is semi-uniform Feller;
\item[{\rm(b)}] the stochastic kernel $\Psi$ on $\S_1\times\S_2$ given $\S_3$ is WTV-continuous;
\item[{\rm(c)}] if $s_3^{(n)}$ converges to $s_3$ in $\S_3,$ then for each closed set $C$ in $\S_1$
\begin{equation}\label{eq:wtvscConv}
\lim_{n\to\infty} \sup_{B\in \B(\S_2)} \left( \Psi(C \times B|s_3^{(n)})-\Psi(C \times B|s_3)\right)= 0;
\end{equation}
\item[{\rm(d)}] if $s_3^{(n)}$ converges to $s_3$ in $\S_3,$ then, for each $A\in\B(\S_1)$ such that $\Psi(\partial A,\S_2|s_3)=0,$
\begin{equation}\label{eq:equivWTV2}
\lim_{n\to\infty} \sup_{B\in \B(\S_2)} | \Psi(A \times B|s_3^{(n)})-\Psi(A \times B|s_3)|= 0;
\end{equation}
\item[{\rm(e)}] if $s_3^{(n)}$ converges to $s_3$ in $\S_3,$ then, for each nonnegative bounded lower semi-continuous function $f$ on $\S_1,$
\begin{equation}\label{eq:equivWTV4}
\ilim_{n\to\infty} \inf_{B\in \B(\S_2)} \left( \int_{\S_1} f(s_1) \Psi(ds_1,B|s_3^{(n)})-\int_{\S_1} f(s_1) \Psi(ds_1,B|s_3)\right)= 0;
\end{equation}
\end{itemize}
{\color{black} and each of these conditions implies continuity in total variation of the marginal kernel $\Psi(\S_1,\,\cdot\,|\,\cdot\,)$ on $\S_2$ given $\S_3.$ }
\end{theorem}


Note that, since $\emptyset\in\B(\S_2),$ \eqref{eq:wtvscConv} holds if and only if
\begin{equation}\label{eq:equivWTV1}
\slim_{n\to\infty} \sup_{B\in \B(\S_2)\setminus\{\emptyset\}} \left( \Psi(C \times B|s_3^{(n)})-\Psi(C \times B|s_3)\right)\le 0,
\end{equation}
and similar remarks are applicable to \eqref{eq:equivWTV2} and \eqref{eq:equivWTV4} with the inequality ``$\ge$'' taking place in \eqref{eq:equivWTV4}.
%
%
Now let $\S_4$ be a Borel subset of a Polish space, and let $\Xi$ be a stochastic kernel on $\S_1\times\S_2$ given $\S_3\times\S_4.$ Consider the stochastic kernel $\Psii$ on $\S_1\times\S_2$ given $\P(\S_3)\times\S_4$ defined by
\begin{equation}\label{eq:extra1}
\Psii(A\times B|\mu,s_4):=\int_{\S_3}\Xi(A\times B |s_3,s_4)\mu(ds_3),
\end{equation}
$A\in\B(\S_1),\,B\in\B(\S_2),\,\mu\in\P(\S_3),\,s_4\in\S_4.$

{\color{black}Note that $\Xi$ is the integrand for $\Psii,$ which justifies the notation $\Psii.$}
The following theorem establishes the preservation of {\color{black}semi-uniform Feller continuity under} the integration  operation in \eqref{eq:extra1}.

\begin{theorem}\label{th:extra}(\cite{FKZJOT})
{\color{black}A} stochastic kernel $\Psii$ on
$\S_1\times\S_2$ given $\P(\S_3)\times\S_4$ is {\color{black}semi-uniform Feller}
if and only if \, $\Xi$ on $\S_1\times\S_2$ given $\S_3\times\S_4$ is {\color{black}semi-uniform Feller}.
\end{theorem}

Let us consider the following assumption.
\begin{Assumption}\label{AssKern}(\cite{FKZJOT})
Let for each $s_3\in\S_3$ the topology on $\S_1$ have a countable base
$\tau_b^{s_3}(\S_1)$ such that
\begin{itemize}
 \item[{\rm(i)}] $\S_1\in\tau_b^{s_3}(\S_1);$
 \item[{\rm(ii)}]  for
each finite intersection $\oo=\cap_{i=1}^ k {\oo}_{i},$ $k=1,2,\ldots,$ of sets
$\oo_{i}\in\tau_b^{s_3}(\S_1),$ $i=1,2,\ldots,k,$
the set of
functions $\fff^\Psi_\oo,$ defined in \eqref{eq:familyoffunctions} with $A=\oo$, is equicontinuous at $s_3.$
\end{itemize}
\end{Assumption}
Let $\S_1,\S_2,$ and $\S_3$ be Borel subsets of Polish spaces, and $\Psi$ on $\S_1\times\S_2$ given $\S_3$ be a stochastic kernel. By Bertsekas and Shreve \cite[Proposition~7.27]{BS}, there exists a stochastic kernel $\Phi$ on $\S_1$ given
$\S_2\times\S_3$ such that
\begin{equation}\label{eq:CBR1}
\Psi(A\times B|s_3)=\int_B\Phi(A|s_2,s_3)\Psi(\S_1,ds_2|s_3),\quad A\in \mathcal{B}(\S_1),\  B\in \mathcal{B}(\S_2),\ s_3\in\S_3.
\end{equation}

The stochastic kernel $\Phi(\,\cdot\,|s_2,s_3)$ on $\S_1$ given
$\S_2\times\S_3$ defines a measurable
mapping $\Phi:\,\S_2\times\S_3 \to\P(\S_1),$ where
$\Phi(s_2,s_3)(\,\cdot\,)=\Phi(\,\cdot\,|s_2,s_3).$ According to Bertsekas and Shreve \cite[Corollary~7.27.1]{BS}, for each $s_3\in
\S_3$ the mapping $\Phi(\,\cdot\,,s_3):\S_2\to\P(\S_1)$ is defined
$\Psi(\S_1,\,\cdot\,|s_3)$-almost surely uniquely in $s_2\in\S_2.$
Let us consider the stochastic kernel  $\phi$ defined by
\begin{equation}\label{eq:CBR2}
\phi(D\times B|s_3):=\int_{B}\h\{\Phi(s_2,s_3)\in D\}\Psi(\S_1,ds_2|s_3),
\end{equation}
$D\in \mathcal{B}(\P(\S_1)),$ $B\in\B(\S_2),$ $s_3\in\S_3,$ where a particular choice of
a stochastic kernel $\Phi$ satisfying (\ref{eq:CBR1}) does not effect the
definition of $\phi$ in (\ref{eq:CBR2}).

In models for decision making with incomplete information, $\phi$ is the transition probability to the set of pairs $(z,y),$  where $z$ is a are posterior probability distribution of a state, and  $y$ is an observation; \eqref{3.7}. Continuity properties of $\phi$ play the fundamental role in the studies of models with incomplete information. Theorem~\ref{th:CBRmain} characterizes such properties, and this is the reason for the title of this section. Let us consider the following assumption.
\begin{Assumption}\label{Ass:H} (\cite{FKZ}) {For a stochastic kernel $\Psi$ on $\S_1\times\S_2$ given $\S_3,$ there} exists a stochastic
kernel $\Phi$ on $\S_1$ given $\S_2\times\S_3$ satisfying
(\ref{eq:CBR1}) such that, if a sequence
$\{s_3^{\left(n\right)}\}_{n=1,2,\ldots}\subset\S_3$ converges to
$s_3\in\S_3$ as $n\to\infty,$ then there exists a
subsequence $\{s_3^{\left(n_k\right)}\}_{k=1,2,\ldots}\subset
\{s_3^{\left(n\right)}\}_{n=1,2,\ldots}$ and a measurable subset $B$ of
$\,\S_2$ such that
\begin{equation}\label{eq:CBR3}
\Psi(\S_1,B|s_3)=1\ \mbox{and}\ \Phi(s_2,s_3^{\left(n_k\right)})\mbox{ converges weakly to }\Phi(s_2,s_3)\ \mbox{for all }s_2\in B.
\end{equation}
In other words, the convergence in \eqref{eq:CBR3} holds $\Psi(\S_1,ds_2|s_3)$-almost
surely.
\end{Assumption}
%

\begin{theorem}\label{th:CBRmain}(\cite{FKZJOT,FKZSIAM})
For a stochastic kernel $\Psi$ on $\S_1\times\S_2$ given $\S_3$ the following conditions are equivalent:
\begin{itemize}
\item[{\rm(a)}] the stochastic kernel $\Psi$ on $\S_1\times\S_2$ given $\S_3$ is semi-uniform Feller;
\item[\rm(b)] the stochastic kernel $\Psi$ on $\S_1\times\S_2$ given $\S_3$  satisfies Assumption~\ref{AssKern};
\item[{\rm(c)}] the marginal kernel $\Psi(\S_1,\,\cdot\,|\,\cdot\,)$ on $\S_2$ given $\S_3$ is continuous in total variation and Assumption~\ref{Ass:H} holds;
\item[{\rm({d})}] the stochastic kernel $\phi$ on $\P(\S_1)\times \S_2$ given $\S_3$ is semi-uniform {Feller.} 
\end{itemize}
\end{theorem}

For a metric space $\S,$ we say that a subset $\F(\S)$ of the set of bounded continuous functions $f:\S\to\R$ \textit{determines weak convergence} on $\P(\S)$  if a sequence of probability measures $\{\mu^{(n)}\}_{n=1,2,\ldots}$ from $\P(\S)$ converges weakly to $\mu\in\P(\S)$ if and only if \eqref{eq:wcEF} holds for all $f\in\F(\S).$ According to \cite[Theorem 6.6, p. 47]{Part}, if a metric space $\S$ is separable, then there exists a countable set $\F(\S)$ of uniformly bounded continuous functions on $\S,$ which determines weak convergence on $\P(\S).$ If a bounded continuous function is added to $\F(\S),$ then the new set also determines weak convergence. Therefore, without loss of generality, we can assume that the function ${\bf I}_\S$ belongs to $\F(\S),$ where ${\bf I}_\S(s)=1$ for all $s\in\S.$  The following assumption is motivated by   \cite[Assumption (M)]{Saldi}; see Section \ref{s4} below for details.
\begin{Assumption}\label{Ass:M}  For a stochastic kernel $\Psi$ on $\S_1\times\S_2$ given $\S_3,$ there exists a countable subset $\F(\S_1)$  of the set of bounded continuous functions $f:\S_1\to\R$ determining weak convergence on $\P(\S_1)$ such that ${\bf I}_{\S_1}\in\F({\color{blue}\S_1}),$ and equality~\eqref{eq:equivWTV3} holds for all $f\in\F(\S_1).$
\end{Assumption}
The following theorem is the central result in this paper.
\begin{theorem}\label{th:CBRLmain}
A stochastic kernel $\Psi$ on $\S_1\times\S_2$ given $\S_3$  is semi-uniform Feller if and only if Assumption~\ref{Ass:M} holds.
\end{theorem}
\begin{proof}
A semi-uniform Feller kernel $\Psi$ on $\S_1\times\S_2$ given $\S_3$ satisfies equality  \eqref{eq:equivWTV3} for all bounded continuous functions $f$ on $\S_1,$ and therefore $\Psi$ satisfy Assumption~\ref{Ass:M}.

Now, let   Assumption~\ref{Ass:M} holds.
%
%
  The assumption ${\bf I}_{\S_1}\in\F(\S_1)$ means that the marginal kernel $\Psi(\S_1, \,\cdot\,|\,\cdot\,)$ is continuous in total variation. Let us fix an arbitrary $s_3\in\S_3.$ Let $f\in\F(\S_1).$ Since the function $f$ is bounded and the marginal kernel  $\Psi(\S_1, \,\cdot\,|\,\cdot\,)$ is continuous in total variation,  \eqref{eq:equivWTV3} and \eqref{eq:CBR1} imply
\begin{equation}\label{eq:neww01}
\begin{aligned}
\lim_{n\to\infty} \sup_{B\in \B(\S_2)} | \int_{B}\int_{\S_1} f(s_1)\Phi(ds_1|s_2,s_3^{\left(n\right)})\Psi(\S_1,ds_2|s_3) \\  -\int_{B}\int_{\S_1} f(s_1)\Phi(ds_1|s_2,s_3)\Psi(\S_1,ds_2|s_3)|= 0
\end{aligned}
\end{equation}
because the family of Borel measurable functions\\ $\{s_2\mapsto\int_{\S_1} f(s_1)\Phi(ds_1|s_2,s_3^{\left(n\right)})\,:\,n=1,2,\ldots\}$ is uniformly bounded on $\S_2$ by the same constant as $f$ is bounded on $\S_1.$ This is equivalent to\\ $\int_{\S_1} f(s_1)\Phi(ds_1|\,\cdot\,,s_3^{\left(n\right)})\to \int_{\S_1} f(s_1)\Phi(ds_1|\,\cdot\,,s_3)$ in $L_1(\S_2,\B(\S_2), \Psi(\S_1,\,\cdot\,|s_3)).$ {Therefore,}
\begin{equation}\label{eqwcef2022}
\int_{\S_1} f(s_1)\Phi(ds_1|\,\cdot\,,s_3^{\left(n_k\right)})\to \int_{\S_1} f(s_1)\Phi(ds_1|\,\cdot\,,s_3) \quad \Psi(\S_1,\,\cdot\,|s_3)\mbox{-a.s., as }k\to \infty,
\end{equation}
for some subsequence $\{n_k\}_{k=1,2,\ldots}$ ($n_k\uparrow\infty$ as $k\to\infty$). Since  \eqref{eqwcef2022} holds for all $f\in\F(\S_1),$  it holds  for all bounded continuous functions $f:\S_1\to\R.$  Thus, Assumption~\ref{Ass:H} holds.  In view of  Theorem~\ref{th:CBRmain}(c), the stochastic kernel $\Psi$ is semi-uniform Feller.
\end{proof}

\section{Semi-Uniform Feller Continuity of Transition Probabilities for MDPCIs}\label{s3}
We start with  the description of the well-known reduction of an MDPII $(\W\times\Y,\A,P,c)$ to an MDPCI (\cite{BS, DY, FKZSIAM, Rh, Yu}). For epoch $t=0,1,\ldots$ consider
the joint conditional probability $R(dw_{t+1}dy_{t+1}|z_t,y_t,a_t)$ on next state $(w_{t+1},y_{t+1})$ given the current posterior state distribution $z_t\in \P(\W),$ observation $y_t\in \Y, $ and the current control action $a_t$ defined by
\begin{equation}\label{3.3}
R(B\times C|z,y,a):=\int_{\W}P(B\times C|w,y,a)z(dw),
\end{equation}
where $B\in \B(\W),$ $C\in\B(\Y),$ $(z,y,a)\in\P(\W)\times\Y\times\A.$ In view of \eqref{eq:CBR1}, there exists a stochastic kernel $H(z,y,a,y')[\,\cdot\,]=H(\,\cdot\,|z,y,a,y')$ on $\W$ given
$\P(\W)\times\Y\times \A\times\Y$ such that
\begin{equation}\label{3.4}
R(B\times C|z,y,a)=\int_CH(B|z,y,a,y')R(\W,dy'|z,y,a),
\end{equation}
where $B\in \B(\W),$ $C\in\B(\Y),$ $(z,y,a)\in\P(\W)\times\Y\times\A.$ The stochastic kernel $H(\,\cdot\,|z,y,a,y')$ introduced in \eqref{3.4} defines a
measurable mapping $H:\,\P(\W)\times\Y\times \A\times \Y {\to}\P(\W).$ Moreover, the mapping $y'\mapsto H(z,y,a,y')$ is defined $R(\W,\,\cdot\,|z,y,a)$-a.s. uniquely for each triplet $(z,y,a)\in \P(\W)\times\Y\times\A.$

Let $\h B$ denotes the \textit{indicator of an event} $B.$ The MDPCI is defined as an MDP with parameters
$(\P(\W)\times\Y,\A,q),$ where
\begin{itemize}
\item[(i)] $\P(\W)\times\Y$ is the state space;
\item[(ii)] $\A$ is the
action set available at all state $(z,y)\in\P(\W)\times\Y;$
 \item[(iii)] $q$ on $\P(\W)\times\Y$
given $\P(\W)\times\Y\times \A$ is a stochastic kernel defined by \eqref{eq:CBR2} with $\S_1:=\W,$ $\S_2:=\Y,$ and $\S_3:=\P(\W)\times\Y\times\A,$  which determines the
distribution of the new
state.  That is,
for $(z,y,a)\in \P(\W)\times\Y\times\A$
and for $D\in \B(\P(\W))$ and $C\in \B(\Y),$
\begin{equation}\label{3.7}
q(D\times C|z,y,a):=\int_{C}\h\{H(z,y,a,y')\in D\}
R(\W,dy'|z, y, a).
\end{equation}

\end{itemize}
 Note that a particular  measurable
choice of a stochastic kernel $H$ from (\ref{3.4}) does not { effect}
the definition of $q$ in (\ref{3.7}).

The transition probability $q,$ which is a stochastic kernel on $\P(\W)\times \Y$ given $\P(\W)\times\Y\times\A,$ defines transition probabilities for MDPCI, and we are interested in establishing its continuity properties.  To do this, it is also useful to write the formula
\begin{equation}\label{3.4m}
P(B\times C|w,y,a)=\int_C H(B|w,y,a,y')P(\W,dy'|w,y,a)
\end{equation}
for $B\in \B(\W),$ $C\in\B(\Y),$ $(w,y,a)\in W\times\Y\times\A,$  which is similar to \eqref{3.4}, and we use the same notation $H$ for the transition probability as in \eqref{3.4} because
$H(B|w,y,a,y')=H(B|\delta_w,y,a,y')$ for all $(w,y,a)\in\W\times\Y\times\A$ almost surely in $P(\W,dy'|w,y,a),$ where $\delta_w$ is the Dirac measure on $\W$ concentrated at $w\in\W.$

In view of Theorem~\ref{th:extra}, the stochastic kernel $P$ is semi-uniform Feller if and only if the stochastic kernel $R$ is semi-uniform Feller.  In view of Theorem~\ref{th:CBRmain}(a,d), the stochastic kernel $R$ is semi-uniform Feller if and only if the stochastic kernel $q$ is semi-uniform Feller.  This leads us to the following theorem.

\begin{theorem}\label{th:mainMDPIINew}(\cite[Theorem 6.2]{FKZSIAM})
Let $(\W\times\Y,\A,P,c)$ be an MDPII, and $(\P(\W)\times\Y,\A,q,\c)$ be its MDPCI. Then the following conditions are equivalent:
\begin{itemize}
\item[{\rm(a)}] the stochastic kernel $P$ on $\W\times\Y$ given
$\W\times\Y\times\A$ is semi-uniform Feller;
\item[{\rm(b)}] the stochastic kernel {$R$ on $\W\times\Y$ given $\P(\W)\times\Y\times \A$ is  semi-uniform Feller;}
\item[{\rm({c})}] the stochastic kernel $q$ on $\P(\W)\times\Y$
given $\P(\W)\times\Y\times \A$ is semi-uniform Feller.
\end{itemize}
\end{theorem}
The most significant fact in Theorem~\ref{th:mainMDPIINew} is that semi-uniform Feller continuity of $P$ is necessary and sufficient for semi-uniform Feller continuity of $q.$
Theorems~\ref{th:equivWTV}-\ref{th:CBRLmain} provide necessary and sufficient conditions for semi-uniform Feller continuity.  Theorem~\ref{th:equivWTV} provides conditions based on the definition of semi-uniform Feller continuity. Theorem~\ref{th:extra} claims preservation of semi-uniform Feller continuity under integration.  In particular, Theorem~\ref{th:extra} implies statement (b) in Theorem~\ref{th:mainMDPIINew}.  Theorems~\ref{th:CBRmain} and \ref{th:CBRLmain} prove that each of the Assumptions~\ref{AssKern} and \ref{Ass:M} is necessary and sufficient for semi-uniform continuity of a kernel.  Theorem~\ref{th:CBRmain} also claims that Assumption~\ref{Ass:H} and the assumption that  the marginal kernel $\Psi(\S_1,\,\cdot\,|\,\cdot\,)$ on $\S_2$ given $\S_3$ is continuous in total variation taken together are necessary and sufficient for semi-uniform Feller continuity of $\Psi.$  Assumption \ref{AssKern} deals with equicontinuity properties of stochastic kernels $\Psi$ considered at certain sets,\ Assumption~\ref{Ass:H} deals with weak continuity of stochastic kernels $\Phi,$ and Assumption~\ref{Ass:M} deals with equicontinuity of integrals for a countable set of functions determining weak convergence.
\section{Continuity of Transition Probabilities for Belief-MDPs}\label{s4}
We recall that Platzman's model is an MDPII whose transition probability $P$ is a stochastic kernel on $\W\times\Y$ given $\W\times\A.$ For Platzman's models and, in particular, for POMDPs, it is possible to consider a completely observable MDP, called belief-MDP, whose state space is $\P(\W),$ and the set of actions is $\A.$   The transition probability $\hat{q}$ for the belief-MDPs is
\begin{equation}\label{eq:filtpef}
\hat{q}(D|z,a):=q(D,\Y|z,a)=\int_{\Y}\h\{H(z,a,y')\in D\}R(\W,dy'|z, a), 
  \end{equation}
where $D\in \B(\P(\W)),$ $z\in\P(\W),$ $a\in\A,$ and $y'\in\Y.$ We recall that for Platzmain's models, including POMDPs, transition probabilities $P$ do not depend on observations $y,$ that is, $P(\cdot,\cdot|w,y,a)=P(\cdot,\cdot|w,a),$ and formulae \eqref{3.3} and \eqref{3.4} become
\begin{equation}\label{3.3ef}
R(B\times C|z,a):=\int_{\W}P(B\times C|w,a)z(dw),
\end{equation}
and
\begin{equation}\label{3.4b}
R(B\times C|z,a)=\int_CH(B|z,a,y')R(\W,dy'|z,a).
\end{equation}
Semi-uniform Feller continuity of the transition probability $q$ implies its weak continuity, which implies weak continuity of its marginal probability $\hat{q}.$ Therefore, the results of Sections~\ref{s2} and \ref{s3} provide sufficient conditions for weak continuity of $\hat{q}.$  In view of Theorem~\ref{th:mainMDPIINew}, semi-uniform Feller continuity of the stochastic kernel $P$ implies weak continuity of $\hat{q}.$

Formula \eqref{3.4m} can be simplified for Platzman's models to
\begin{equation}\label{3.4ma}
P(B\times C|w,a)=\int_C H(B|w,a,y')P(\W,dy'|w,a),
\end{equation}
where formula \eqref{3.4ma} is related to formula \eqref{3.4b} in the same way  \eqref{3.4m} is related to \eqref{3.4}. In particular, the relation between the kernel $H$  on $\W\times\Y$  given $\W\times\Y$ in \eqref{3.4ma} and the kernel $H$ on $\W\times\Y$ given $\P(\W)\times\Y$ in \eqref{3.4b}   is $H(B|\delta_w,a,y')=H(B|w,a,y')$ for all $(w,a)\in\W\times\A$ almost surely in $P(\W,dy'|w,a).$

According to Theorem~\ref{th:mainMDPIINew}, there are two  approaches to prove semi-uniform Feller continuity of the kernel $q$: (i) prove semi-uniform continuity of
 $P,$ and (ii) prove semi-uniform continuity of $R.$ The kernel $R$ defines the kernel $\hat{q}$ via \eqref{eq:filtpef}, and kernel
$R$ was used to prove weak continuity of $\hat{q}$ in several references including \cite{FKZ, HL, Saldi}. However, it is typically easier to use approach (i) than (ii)  to prove  semi-uniform Feller continuity of $q.$ In particular, formula \eqref{3.4ma} is useful for verifying Assumption~~\ref{Ass:H} for the kernel $P.$ 
%
%

In the literature on POMDPs, the transition probability $\hat{q}$ is usually defined by the right-hand side of \eqref{eq:filtpef}, and the transition probability $q$ is not considered.  Here and in \cite{FKZSIAM} we consider $q$ because its weak continuity implies weak continuity of $\hat{q}.$ The transition probability $q$ is important for MDPCIs. Platzman's models including POMDPs are particular cases of MDPIIs, and MDPCIs can be also constructed for them.  The state space of an MDPCI is $\P(\W)\times\Y.$ However, if one-step costs do not depend on observations, neither transition probability between belief states $z\in\P(\W)$ nor  costs depend on observations $y\in\Y.$ For such problems, the set $\Y$ contains non-essential information, and, therefore, it is sufficient to consider only the state space $\P(\W)$ for belief-MDPs for Platzman's models including POMDPs when costs do not depend on observations; see \cite{FKZSIAM} for details.  The general theory for such reductions is described in \cite{Fe}. The original development of that theory was motivated by Continuous-Time Markov Decision Processes \cite{PZ} and their reduction to discrete time \cite{Fe12}.

Recall that ${\rm POMDP}_1$ is Platzman's model with the transition probability
\[
P(B\times C|w,a)=P_1(B|w,a)Q_1(C|w,a),
\]
where $B\in \B(\W),$ $C\in\B(\Y),$ $w\in\W,$ $a\in\A,$ $P_1$ is a stochastic kernel on $\W$ given $\W\times \A,$ and $Q_1$ is a stochastic kernel on $\Y$ given $\W\times \A.$  Thus, $P_1$ is the transition probability for the  MDP with hidden states, and $Q_1$ is the observation probability.
 For a ${\rm POMDP}_1$ semi-uniform Feller continuity of $P$ is equivalent to the validity of the following properties: the transition probability $P_1$ is weakly continuous, and the observation probability $Q_1$ is continuous in total variation \cite[Corollary 6.10]{FKZSIAM}.

Recall that ${\rm POMDP}_2$ is Platzman's model with the transition probability
\begin{equation}\label{eqtrbmdp2}
P(B\times C|w,a):=\int_B Q_2(C|a,w')P_2(dw'|w,a),
\end{equation}
where $B\in \B(\W),$ $C\in\B(\Y),$ $w\in\W,$ $a\in\A,$ $P_1$ is a stochastic kernel on $\W$ given $\W\times \A,$ and $Q_2$ is a stochastic kernel on $\Y$ given $\A\times \W.$  Thus, $P_2$ is the transition probability for the  MDP with hidden states, and $Q_2$ is the observation probability.

For ${\rm POMDP}_2$ semi-uniform Feller continuity of $P$ holds in the following two cases \cite[Corollary 6.10]{FKZSIAM}:
\begin{itemize}
\item[{\rm(i)}] the transition probability $P_2$ is weakly continuous, and the observation probability $Q_2$ is continuous in total variation;
\item[{\rm(ii)}] the transition probability $P_2$ is continuous in total variation, and the observation probability $Q_2(\,\cdot\,|a,\,\cdot\,)$ is continuous in total variation in the control parameter $a\in\A.$
\end{itemize}
Thus, if the transition probability $P_i$ is weakly continuous, and the observation probability $Q_i$ is continuous in total variation, then the transition probability $\hat{q}$ is weakly continuity for ${\rm POMDP}_i,$ $i=1,2.$ In addition, if the transition probability $P_2$ is continuous in total variation, and observation probability $Q_2(\,\cdot\,|a,\,\cdot\,)$ is continuous in total variation in the control parameter $a,$ then the transition probability $\hat{q}$ is weakly continuity for  ${\rm POMDP}_2.$

Sufficiency of condition (i) for weak continuity of the transition kernel $\hat{q}$ for a ${\rm POMDP}_2$ was proved directly in \cite{FKZ}. Another proof of this fact was provided in \cite{Saldi}, where also the following sufficient condition,  for weak continuity of $\hat{q}$ was established:
 \begin{itemize}
\item[{\rm(iii)}] the transition probability $P_2$ is continuous in total variation, and the observation probability $Q_2$ does not depend on the control parameter $a.$
\end{itemize}
Condition (ii) is a generalization of condition (iii).

Thus, for ${\rm POMDP}_1$ weak continuity of $P_1$ and continuity of $Q_1$ in total variation are the necessary and sufficient conditions for semi-uniform Feller continuity of $P.$  For ${\rm POMDP}_1$ statements (i) and (ii) provide sufficient conditions for weak continuity of $P$.  The natural question is whether  conditions (i) and (ii) taken together are necessary?  Example~\ref{exa:1} provides the negative answer to this question.  Therefore, criteria for semi-uniform Feller continuity are important for studying ${\rm POMDP}_2.$

Let us consider an example of  ${\rm POMDP}_2$ with a semi-uniform Feller continuous kernel $P$ which falls neither into case (i) nor into case (ii).
\begin{example}\label{exa:1}
The transition kernel $P_2$ on $\W$ given $\W\times\A$ is weakly continuous, but it is not continuous in total variation, the observation kernel $Q_2$ on $\Y$ given $\A\times\W$ does not depend on the control parameter $a,$ and it is not continuous in total variation, and the transition kernel $P$ on  $\W\times\Y$ given $\W\times\A$ is semi-uniform Feller continuous.

{\rm Let $d_+:=\max\{d,0\},$ and $d_-:=\min\{d,0\}$ for each $d\in\mathbb{R}.$ We set $\W=\Y=\A:=\mathbb{R},$ $P_2(B|w,a):={\bf I}\{w_+\in B\},$ and $Q_2(C|w):={\bf I}\{w_-\in C\},$ $w,a\in \mathbb{R},$ $B,C\in\mathcal{B}(\mathbb{R}).$ Then $\int_\W f(w')P_2(dw'|w,a)=f(w_+)$ and $\int_\Y g(y)Q_2(dy|w)=g(w_-)$ for bounded continuous functions $f$ and $g.$   Stochastic kernels $P_2$ and $Q_2$ are obviously weakly continuous at each $w\in\mathbb{R},$ but each of them is not continuous in total variation at $w=0.$ Moreover, direct calculations imply that $P(B\times C|w,a)={\bf I}\{w_+\in B\}{\bf I}\{0\in C\},$ $B,C\in\mathcal{B}(\mathbb{R}),$ $w,a\in\mathbb{R},$ and $P$ is semi-uniform Feller continuous because
 for each sequence $\{w^{(n)}\}_{n=1,2,\ldots}\subset\R$ that converges to $w\in\R$ and for each bounded continuous function $f$ on $\R,$
\[
\begin{aligned}
\lim_{n\to\infty}& \sup_{C\in \B(\R)} \left| \int_{\R} f(w') P(dw',C|w^{(n)})-\int_{\R} f(w') P(dw',C|w)\right|\\
&= \lim_{n\to\infty} \sup_{C\in \B(\R)} {\bf I}\{ 0 \in C \} \left|  f(w_+^{(n)})-f(w_+)\right|= 0,
\end{aligned}
\]
where the last equality follows from continuity of $f$ on $\R.$
}
\end{example}

\begin{Remark}
{\rm Since $q$ is semi-uniform Feller if and only if $P$ is semi-uniform Feller, then $q$ and $\hat{q}$ are weakly continuous if $P$ is semi-uniform Feller.
 However,
it is possible that $\hat{q}$ is weakly continuous, but $P$ is not semi-uniform Feller.  For example, let us present an MDP with the state space $\W,$ and the action space $\A,$ and transition probability $p(B|w,a)={\bf I}\{w\in B\}$  as ${\rm POMDP}_2$  with $\Y=\W,$ $P_2(B|w,a)={\bf I}\{w\in B\},$  and $Q_2(C|a,w)={\bf I}\{w\in C\},$ where  $w\in\W,$  $a\in\A,$  $y\in\Y,$ $B\in\B(\W),$ and $C\in\B(\Y).$  Then $P(B\times C|w,a)={\bf I}\{w\in B\cap C \},$ and the kernel $P$ is not semi-uniform Feller.  It is easy to see that $\hat{q}$ is weakly continuous in this example.  In particular,   for this example $H(B|z,a,y)= {\bf I}\{y\in B\}$ satisfies \eqref{3.4b}.  The kernel $H$ is weakly continuous, and together with weak continuity of $P_2$ and $Q_2$ this is a sufficient condition for weak continuity of $\hat{q},$ see e.g.,  \cite[p. 90]{HL} or \cite[Theorem 3.2]{FKZ}.
}
\end{Remark}

Assumption~\ref{AssKern} was introduced in \cite{FKZJOT}, and its stronger version, when the base $\tau_b^{s_3}(\S_1)$ does not depend on $s_3,$ was introduced in \cite{Steklov} to study MDPIIs.  Assumption~\ref{Ass:H} was introduced in \cite{FKZ} for the transition probability $R$ defined in \eqref{3.3} for the transition probability $P$ defined in  \eqref{eqtrbmdp2}. Assumption~\ref{Ass:M} is relevant to Assumption (M) introduced in \cite{Saldi} for the transition probability $R$  as an alternative to Assumption \ref{Ass:H} for a sufficient condition of weak continuity of the transition probability $\hat{q}$ for ${\rm POMDP}_2.$
In terms of this paper, Assumption (M) from \cite{Saldi} can be formulated in the following form.

\noindent{\bf Assumption (M).} (\cite{Saldi}) For a countable set $\F(\W)=\{f_m\}_{m\ge 1}$ of uniformly bounded continuous functions $f:\W\to \R$ such that:

\begin{itemize}\item[(a)] ${\bf I}_{\W}\in\F;$  \item[(b)] $\F(\W)$ metrizes the weak topology on $\P(\W)$ with the metric
\begin{equation}\label{eq:mtrho}
\rho(\mu,\nu):=\sum_{m=1}^\infty 2^{-m}\left|\int_{\W} f_m(w)\mu(dw)-\int_{\W} f_m(w)\nu(dw) \right|,
\end{equation}
\item[(c)] equicontinuity property \eqref{eq:equivWTV3} holds for all $f\in\F(\S_1)$ with $s_1=w,$ $\S_1=\W,$ $\S_2=\Y,$ $\S_3=\P(\W)\times\A,$ and $\Psi=R.$
\end{itemize}

Assumptions (M) can be viewed as an  implementation of Assumption~\ref{Ass:M} for particular spaces.  The following two differences are not essential:
\begin{itemize}
\item[] Assumption~\ref{Ass:M} states that the functions in $\F(\S_1)$ are bounded, and Assumption (M) assumes that the functions in $\F(\W)$ are uniformly bounded;
\item[] Assumption~\ref{Ass:M} states that the set of functions $\F(\S_1)$ determines the topology of weak convergence, while Assumption~(M) states the metric $\rho$ defined in \eqref{eq:mtrho} metrizes the topology of weak convergence  on $\P(\W).$
\end{itemize}
Indeed, the family $\F(\S_1)=\{f_m\}_{m\ge 1}$ in Assumption~\ref{Ass:M} consists of bounded functions.  This means that $\sup_{s_1\in\S_1}|f_m(s_1)|\le L_m<+\infty$ for all $m=1,2,\dots.$  Then $\{f_m/\max\{L_m,1\}\}_{m=1,2,\ldots}$ is the set of uniformly bounded functions satisfying all the conditions in Assumption~\ref{Ass:M}.  In addition, when $\S_1=\W,$ the condition that the set $\F(\S_1)$ determines weak convergence on $\P(\S_1)$ and the condition that the metric $\rho$ defined in \eqref{eq:mtrho}  metrizes the topology of weak convergence  on $\P(\W)$ are obviously equivalent since $\W$ is a metric space.

It was observed in \cite{Saldi} for that ${\rm POMDP}_2$  that Assumption (M)
is more general than assumptions~(i) and~(iii) stated in this section.  Indeed, as follows from \cite[Corollary 6.10]{FKZSIAM} and Theorems~\ref{th:CBRLmain},~\ref{th:mainMDPIINew},  assumptions~(i)--(iii) from this section are sufficient conditions for semi-uniform Feller continuity of each of the transition probabilities $P,$ $R,$ and $q,$  while Assumption (M) is the necessary and sufficient conditions  for semi-uniform Feller continuity of $P,$ $R,$ and $q.$

\bibliographystyle{elsarticle-num}

\begin{thebibliography}{9}
%


\bibitem{Ao} Aoki, M. (1965) Optimal control of partially observable Markovian systems. \emph{J. Franklin Inst.} 280(5): 367--386.

\bibitem{As}
 {\AA}str\"om, K.J. (1965). Optimal control of Markov processes with incomplete state information. \emph{Journal of Mathematical Analysis and Applications} 10: 174--205.


 \bibitem{BR} B\"auerle, N., Rieder, U. (2011) \textit{Markov Decision Processes with Applications to Finance,} Springer-Verlag, Berlin.

\bibitem{BS} Bertsekas, D.P.,  Shreve S.E. (1978) \textit{Stochastic Optimal Control: The Discrete-Time
Case,} Academic Press, New York 



\bibitem{Dy}  Dynkin, E.B. (1965) Controlled random sequences. \emph{Theory
Probab. Appl.} 10(1): 1--14.

\bibitem{DY} Dynkin, E.B.,   Yushkevich, A.A. (1979) \textit{{C}ontrolled
{M}arkov {P}rocesses,} Springer-Verlag, New York.

\bibitem{Fe} Feinberg, E.A. (2005) On essential information in sequential decision processes. Math. Meth. Oper. Res. 62, 399--410.

\bibitem{Fe12}
Feinberg E (2012) Reduction of discounted continuous-time mdps with unbounded
  jump and reward rates to discrete-time total-reward {MDP}s. Hernandez D,
  Minjarez A, eds., \emph{Optimization, Control, and Applications of Stochastic
  Systems {\rm (Birkh\"auser/Springer, New York)}}, 201--213.

\bibitem{FKL2} Feinberg, E.A.,  Kasyanov, P.O., Liang, Y. (2020) Fatou's lemma in its classical form and Lebesgue's
convergence theorems for varying measures with applications to Markov decision processes. \emph{Theory Probab. Appl.} 65(2): 270--291.


\bibitem{Steklov} Feinberg, E.A., Kasyanov, P.O., Zgurovsky, M.Z. (2014) Convergence of probability measures and Markov decision models with incomplete information,
\emph{Proceedings of the Steklov Institute of Mathematics,} 287 (1), 96--117.

\bibitem{UFL} Feinberg, E.A.,  Kasyanov, P.O., Zgurovsky, M.Z. (2016) Uniform Fatou's lemma,  \emph{Journal of Mathematical Analysis and Applications}, 444(1), 550--567.

\bibitem{FKZ} Feinberg, E.A., Kasyanov, P.O., Zgurovsky, M.Z. (2016) Partially observable total-cost Markov decision processes with weakly continuous transition probabilities, \emph{Math. Oper. Res.}, 41(2), 656--681.

    \bibitem{FKZJOT} Feinberg, E.A.,  Kasyanov, P.O., Zgurovsky, M.Z. (2021) Semi-uniform
Feller kernels,  	arXiv:2107.02207.

\bibitem{FKZSIAM} Feinberg, E.A.,  Kasyanov, P.O., Zgurovsky, M.Z. (2021) Markov decision processes with incomplete information and semi-uniform Feller transition probabilities,  arXiv:2108.09232, {\color{black} to appear in \emph{SIAM Journal on Control and Optimization}}.


{\color{black}
\bibitem{HL}  Hern\'{a}ndez-Lerma, O. (1989) \textit{Adaptive Markov Control Processes,} Springer-Verlag, New York.
}

\bibitem{Saldi} Kara, A.D., Saldi, N., Y\"uksel, S. (2019) Weak Feller property of non-linear filters. \emph{Systems $\&$ Control Letters} 134: 104512.



\bibitem{Part}  Parthasarathy, K.R. (1967)  \textit{Probability Measures on Metric Spaces,} Academic Press, New York.

\bibitem{PZ}
Piunovskiy A, Zhang Y (2020) \emph{Continuous-{T}ime {M}arkov {D}ecision
  {P}rocesses: {B}orel {S}pace {M}odels and {G}eneral {S}trategies} (Springer
  Nature, Cham, Switzerland).





\bibitem{Plat} Platzman, L.K. (1980) Optimal infinite-horizon undiscounted control of finite probabilistic systems. \emph{SIAM Journal on Control and Optimization} 18(4): 362--380.


\bibitem{Rh}  Rhenius, D. (1974) Incomplete information in Markovian decision models. \emph{Ann. Statist.} 2(6): 1327--1334.


{\color{black}
\bibitem{RS} Runggaldier, W.J.,   Stettner, L. (1994)
\emph{Approximations of Discrete Time Partially Observed Control Problems,}
Applied Mathematics Monographs CNR, Giardini Editori, Pisa.}


{\color{black}\bibitem{sha} Sch\"al, M. (1975) On dynamic programming: compactness of the space of policies. \emph{Stoch. Process.  Appl.} 3:345--364.}

\bibitem{Shi64}{\color{black}
Shiryaev, A.N. (1964) On the theory of decision functions and control by an observation process with incomplete data. \textit{Transactions of the Third Prague
Conference on Information Theory, Statistical Decision Functions,
Random Processes} (Prague, 1962), pp.~657-681 (in Russian); Engl.
transl. in \textit{Select. Transl. Math. Statist. Probab.}
6(1966), 162-188.}


\bibitem{Shi}
Shiryaev, A.N. (1967) Some new results in the theory of controlled
random processes. \textit{Transactions of the Fourth Prague
Conference on Information Theory, Statistical Decision Functions,
Random Processes} (Prague, 1965), pp.~131-201 (in Russian); Engl.
transl. in \textit{Select. Transl. Math. Statist. Probab.}
8(1969), 49-130.


\bibitem{Yu} Yushkevich AA (1976) Reduction of a controlled Markov model with incomplete data to a problem with complete information
in the case of Borel state and control spaces. \emph{Theory
Probab. Appl.} 21(1): 153--158.
\end{thebibliography}

\end{document}